\documentclass[10pt]{article}
\usepackage{amssymb,amsfonts,amsmath,latexsym,verbatim,amscd,amsthm,graphicx,epsfig}
\title{Self-similar expanding solutions \\ for the planar network flow}
\author{Rafe Mazzeo \thanks{Department of Mathematics, Stanford University, 
Stanford, CA 94305; Email: mazzeo@math.stanford.edu; Supported by the NSF under
grant DMS-050579}\\ Stanford University \and Mariel Saez \thanks{Max Planck Institute 
for Gravitational Physics (Albert Einstein Institute), Am M\"uhlenberg 1, D-14476 Golm, 
Germany;  Email: mariel@aei.mpg.de}
\\ Max Planck Institut f\"ur Gravitationsphysik}

\newcommand{\RR}{\mathbb R}
\newcommand{\ZZ}{\mathbb Z}
\newcommand{\e}{\epsilon}

\newcommand{\calF}{\mathcal F}

\newcommand{\calU}{\mathcal U}

\newcommand{\calC}{\mathcal C}

\newcommand{\calP}{\mathcal P}
\newcommand{\del}{\partial}

\newcommand{\SSS}{\mathbb S}

\newtheorem{proposition}{Proposition}[section]
\newtheorem{corollary}{Corollary}[section]
\newtheorem{lemma}{Lemma}[section]
\newtheorem{definition}{Definition}[section]

\newtheorem*{main}{Main Theorem}

\date{April 23, 2007}

\begin{document}
\bibliographystyle{plain}

\maketitle

\begin{abstract}
We prove the existence of self-similar expanding solutions of the curvature flow
on planar networks where the initial configuration is any number of half-lines
meeting at the origin. This generalizes recent work by Schn\"urer and Schulze
which treats the case of three half-lines. There are multiple solutions, and these
are parametrized by combinatorial objects, namely Steiner trees with respect
to a complete negatively curved metric on the unit ball which span $k$ specified points 
on the boundary at infinity. We also provide a sharp formulation of the regularity 
of these solutions at $t=0$. 
\end{abstract}

\section{Introduction}
The detailed analysis of the curve-shortening flow for embedded closed curves in the plane was an early
success in the field of geometric flows. It is not hard to extend this theory to include curves with boundary 
which are either fixed (a Dirichlet condition) or constrained to lie on the boundary of a convex domain,
for example (a Neumann condition). 

Slightly more generally, one might consider the flow by curvature for networks of curves. 
\begin{definition} A planar network is a finite union of embedded arcs and properly embedded half-lines 
$\{\gamma_i\}$ such that for each $i \neq j$, $\gamma_i \cap \gamma_j$ is either empty or else consists
of one (or both) boundary points of each curve. Each boundary points of every $\gamma_i$ is either a boundary or interior vertex. These intersections are called the (interior) vertices
of the network. The boundary of the network consists of the set of points which are endpoints of exactly
one of these curves. The number of curves intersecting at each interior vertex is called the valence of 
that vertex. We always assume that all curves are at least $\calC^2$
\end{definition}
One must formulate the evolution equation for the flow near vertices weakly, and this is commonly done
using Brakke's ideas \cite{Br}. However, certain features of the short-time existence, and many other 
aspects of the long-time existence and convergence, for this network flow have proved elusive. 
The first thorough analytic attack on this problem was undertaken by Mantegazza, Novaga and 
Tortorelli \cite{moybycm} several years ago. Their results are primarily directed at networks which 
consist of only three arcs meeting at equal angles at a single interior vertex. The primary difficulties
include the choice of boundary conditions at this vertex, and later (for more complicated initial
networks) the possibility of creation and annihilation of such vertices during the flow. They obtain
some interesting convergence results under certain hypotheses, but many difficult questions remain.

Let us say that a network is {\it regular} if each interior vertex is trivalent and the curves meet at 
equal angles (of $2\pi/3$) there. The expected behaviour of a network under the Brakke flow which is 
regular at time $t=0$ is that it remains regular at most instants of time. More generally, a vertex is 
stable if the unit tangents to all incoming curves at that vertex sum to zero. Stability of vertices
is also preserved under (some choice of) Brakke flow; however, in this paper we shall focus on
flows where the vertices are generically regular, in the sense above.  

It is interesting to move beyond the restricted class of initial configurations considered in \cite{moybycm}
and consider the flow starting at a more general network with multivalent vertices. One motivation is
that such `nonregular' networks may appear at discrete values of time when vertices collide, so it
is important to understand how to flow past them. This paper takes the first step in proving short-time 
existence for the curvature flow on networks when the initial configuration is one of these more
general networks. We prove here the existence of self-similar solutions, i.e.\ expanding solitons for 
the flow, when the initial network is a finite union of half-lines intersecting at the origin.
The solution with this initial condition is far from unique, but we are able to describe the set of 
all solutions which are regular when $t>0$. Finally, we also present a somewhat new perspective 
which leads to a sharp regularity statement at `irregular vertices' at time $0$. In a forthcoming 
sequel to this paper we shall apply our results here to prove short-time existence for the network 
flow starting from fairly general initial networks. 

Our approach is inspired by a quite recent paper by Schn\"urer and Schulze \cite{selsimos}, in which
they prove the existence and uniqueness of a self-similar solution when the initial condition consists 
of three half-lines meeting at the origin, but not necessarily in equal angles. Their solution is regular 
for $t>0$ and remains a union of three properly embedded arcs meeting at a common vertex. They do not 
state the precise trajectory of this vertex. In contrast, the self-similar solutions here, which start 
from a union of at least four half-lines meeting at $0$, immediately break up into a regular network with 
multiple interior vertices and remain so for all later times. As we explain, the trajectories of these vertices 
are easy to determine. This `explosion' of a nonregular vertex into a more complicated network provides the 
model for the short-time existence for the general network flow, and the multiplicity of self-similar
solutions corresponds to the nonuniqueness of solutions with a given initial condition. 

Imposing self-similarity is tantamount to a dimension reduction of the equation, which transforms 
this problem into an ODE. Schn\"urer and Schulze derive certain convexity properties of solutions 
of this ODE, which were key to their analysis. However, there is a somewhat broader and more natural 
geometric picture which we explain here, that solutions of this ODE are geodesics for a certain complete 
metric on the plane, and the curvature properties of this metric provide a concise 
explanation for those convexity properties. Furthermore, the identification of these curves with 
geodesics allows us to use variational arguments to prove the existence of the more complicated 
regular networks which provide the solutions to our problem. 

To state our main result, let us introduce some notation. Let $B$ be the ball, regarded as the
stereographic compactification of $\RR^2$.  A union of $k$ half-lines $C_0$ meeting at the
origin in $\RR^2$ determines a finite collection of points $p_1, \ldots, p_k \in \del B$.
Define the metric
\[
g = e^{x^2+y^2}(dx^2 + dy^2).
\]
This is complete and negatively curved, with curvature tending to $0$ at infinity. 
We can now state our main result.

\begin{main}
Let $C_0$ be a finite union of half-lines in the plane meeting at $0$, and $p_1, \ldots, p_k$ the 
corresponding points on $\del B$. The set of self-similar solutions of the curve-shortening flow with 
initial condition $C_0$ for which the network is connected at every time $t>0$ is in bijective 
correspondence with the set of possibly disconnected regular networks on $B$, each arc of which 
is a geodesic for $g$, with boundary the $k$ prescribed points at infinity. There always exists 
at least one (and often very many) connected geodesic Steiner tree with these
asymptotic boundary values.  Finally, these self-similar solutions lift to a smooth family of
networks on the parabolic blowup of $\RR^2 \times \RR^+$ at $x=y=t=0$. 
\end{main}

We remark that Steiner trees for the metric $g$ are most likely combinatorially the same as 
Steiner trees for the hyperbolic metric. In other words, given a deformation within the class 
of complete negatively curved metrics between $g$ and the standard hyperbolic metric on $B$,
one should be able to continuously deform each of the regular geodesic networks we have
found here to a regular geodesic network in hyperbolic space without destroying or creating
new interior vertices. Alternately, any direct synthetic procedure which produces Steiner trees 
in ${\mathbb H}^2$ should adapt directly to produce analogous objects for $g$. 
This correspondence would be useful if there were some effective way of enumerating
complete Steiner trees with given boundary in the hyperbolic plane, but this is known
to be a rather difficult problem, so it is probably better to simply note the
resemblance between these two settings and leave it at that.

This general picture was understood qualitatively by Brakke, and hinted at in an appendix to his 
book \cite{Br}. Unbeknownst to us when we were doing this work, many of the specific facts 
presented here were discovered in slightly different forms by Tom Ilmanen and Brian White in 
the mid '90's. Some discussion of this appears in \cite{ACI} and \cite{I}, but the focus 
in those papers is mostly on the higher dimensional case. We hope that this independent and 
more elementary discussion of the one-dimensional case, along the lines of \cite{selsimos}, 
is not unwelcome. There are some interesting new points too, including the enumeration of
self-similar expanding solutions in terms of (nonelementary) combinatorial data, i.e.\ the
number of Steiner trees in $(\RR^2,g)$ spanning the $k$ given points at infinity, and the
formulation of regularity. The first author wishes to thank Brian White for several very
helpful conversations, and in particular for explaining certain aspects of the Brakke flow,
and more importantly, `size minimization' in the class of flat chains mod $k$, which provides 
a shortcut to the existence result of \S 3 which circumvents a more explicit but longer
synthetic approach. The second author wishes to thank Marilyn Daily and Felix Schulze for 
pointing out helpful references.

The next section describes solutions of the dimension reduced equation and the geometry of the metric $g$;
the main existence result for regular networks with prescribed asymptotes is proved in \S 3; finally,
we make some remarks about regularity at $t=0$ in \S 4, explaining the last assertion of the main theorem.

\section{Self-similar solutions of curve-shortening flow}
Let $\gamma_0$ be an immersed curve in the plane. The curve-shortening flow with initial condition $\gamma_0$ 
is the evolution leading to the family of curves $\gamma_t$, $t \geq 0$, where
\[
\frac{d\,}{dt}\gamma_t = \kappa(\gamma_t) \nu(\gamma_t);
\]
here $\kappa$ is the curvature and $\nu$ the unit normal to $\gamma_t$. The well-known theorem of Grayson 
asserts that if $\gamma_0$ is closed and embedded, then $\gamma_t$ remains embedded and 
shrinks to a point at some finite time $T$. Moreover, for an appropriate choice of `center' $P$, 
$(T-t)^{-1/2}(\gamma_t - P)$ converges in $\calC^\infty$ to the circle $S^1$. 

We are being sloppy here and conflating the curves $\gamma_t$ with embeddings $F_t$ from $S^1$ (for example) 
to $\RR^2$ with image $\gamma_t$. 
It is frequently convenient to consider a modified flow equation which includes an extra tangential term; 
the flow leads to the same family of curves, but alters their  parametrizations. We mostly work with 
parametrized curves below, and the equations of motion will have an extra tangential term, i.e.\  have the
form $\dot{\gamma}_t = \kappa \nu + f \gamma'$ for some function $f$. 

A solution $\{\gamma_t\}_{t \geq 0}$ of this equation is called self-similar, or a soliton, if $\gamma_t$ is 
similar to $\gamma_0$ for all $t > 0$. This notion can be defined for curvature flows whenever the ambient
space has a Killing field \cite{HS}. We are particularly interested in solutions for which $\gamma_t$ is 
simply a dilation of $\gamma_0$, $\gamma_t = \lambda(t) \gamma_0$ for some function $\lambda(t)$. An
alternate way to phrase this uses the family of parabolic dilations $D_\lambda$ on $\RR^2 \times \RR^+$, $\lambda > 0$:
\[
D_\lambda(x,y,t) = (\lambda x, \lambda y, \lambda^2 t).
\]
Any solution of the curve-shortening flow determines a `world-sheet'
\[
\Gamma = \bigcup_{t \geq 0} \gamma_t \times \{t\} \subset  \RR^2 \times \RR^+,
\]
and self-similarity is equivalent to the requirement that $D_\lambda(\Gamma) = \Gamma$ for all $\lambda > 0$. 
(This corresponds to the expanding case; shrinking self-similar solutions have world-sheets
$\Gamma \subset \RR^2 \times \RR^-$, so that $t=0$ is the time of extinction, or at least, the boundary
of $\Gamma$.) It is clear that such a 
solution, if it exists, is determined by the set $\gamma_{1/2} = \Gamma \cap \{t=1/2\}$ (the reason
for using $t=1/2$ rather than $t=1$ is to make various other equations neater); furthermore, 
$\gamma_0$ is simply the (unique) tangent cone at infinity of $\gamma_{1/2}$. 

Self-similarity transforms the curve-shortening flow into a stationary equation for this curve $\gamma_{1/2}$, 
which is simply
\begin{equation}
\kappa = (x,y) \cdot \nu.
\label{eq:csre}
\end{equation}
To prove this, fix a parametrization $F(u)$ of (some piece of) the curve at $t=1/2$. Then, 
\[
F_t(u) = \lambda(t)F(u/\lambda(t))
\]
parametrizes the corresponding part of that curve at any other time, where $\lambda(t)$ is
some function to be determined with $\lambda(1/2) = 1$. Setting this expression into the equation yields
\[
\frac{d\,}{dt}F_t(u)\cdot \nu = \left(\dot{\lambda}F(u/\lambda) - \frac{\dot{\lambda}}{\lambda}
\frac{u}{\lambda}F'(u/\lambda)\right) \cdot \nu = \dot{\lambda}F(u/\lambda)\cdot \nu
\]
and
\[
\kappa = \frac{1}{\lambda|F'(u/\lambda)|^2}F''(u/\lambda) \cdot \nu,
\]
hence
\begin{equation}
\frac{1}{|F'(u/\lambda)|^2}F''(u/\lambda) \cdot \nu = \lambda(t)\dot{\lambda}(t) F(u/\lambda)\cdot \nu. 
\label{eq:csre2}
\end{equation}
Here $\dot{\lambda}$ is the derivative with respect to $t$. This can hold for all $u$ and $t$ if and only
if  $\lambda \dot{\lambda} = c$ is constant, so $\lambda^2 = 2ct + c'$. Since $\lambda(0) = 0$ and 
$\lambda(1/2) = 1$, we get $\lambda(t) = \sqrt{2t}$, and so $F_t(u) = \sqrt{2t}F(u/\sqrt{2t})$. 
Setting this in (\ref{eq:csre2}) gives (\ref{eq:csre}).

In any case, we can now reformulate our problem as the
\begin{proposition}
A self-similar solution of the curve-shortening flow with initial condition $C_0$, the union
of a finite number of half-lines meeting at $0$, is equivalent to a regular network of curves
in $\RR^2$, each of which is a solution to (\ref{eq:csre}), with tangent cone at infinity 
equal to $C_0$.
\label{pr:main}
\end{proposition}

In the remainder of this section, we determine all solutions of (\ref{eq:csre}). 

Note first that there is a distinguished subset of solutions, namely the collection of all straight 
lines through the origin in $\RR^2$. Each such line is clearly a solution, since its curvature
is zero and the position vector and tangent vector are always multiples of one another, hence orthogonal
to the normal at each point. This set of lines gives the full set of solutions passing through
the origin, and conversely, any solution passing through $0$ is a straight line. In fact, even
more is true: any solution $\gamma$ which has tangent vector a multiple of the position vector
at any point is one of these straight lines. 

The most convenient parametrization for any other solution is as a normal graph. This is because,
following the last remark above, the tangent never points in the radial direction. Thus we seek
a function $r(\theta)$ so that the curve is the image of the map $F(\theta) = r(\theta) R(\theta)$, where 
\[
R(\theta) = (\cos \theta, \sin \theta), \qquad N(\theta) = (-\sin \theta, \cos \theta).
\]
We calculate
\[
F' = r' R + r N, \qquad 
F''= (r''- r) R + 2r' N,
\]
so in particular
\[
T = \frac{1}{\sigma}F' = \frac{r'}{\sigma}R + \frac{r}{\sigma}N, \qquad \mbox{and} \qquad 
\nu = -\frac{r}{\sigma}R + \frac{r'}{\sigma}N
\]
are the unit tangent and normal; here $\sigma = \sqrt{r^2 + (r')^2}$. Hence 
\[
\kappa \nu = \frac{1}{\sigma}\left(\frac{1}{\sigma}F'\right)' = 
\frac{1}{\sigma^2}F'' - \frac{\sigma'}{\sigma^3}F'.
\]
Rewriting this in terms of $r$, $R$ and $N$, we see that
\[
\kappa \nu  = \left(\frac{r'' - r}{\sigma^2} - \frac{(\sigma')^2}{\sigma^3} r'\right) R +
\left( \frac{2r'}{\sigma^2} - \frac{\sigma'}{\sigma^2}r\right) N,
\]
so finally
\begin{equation}
\kappa = \frac{1}{\sigma^3}\left( 2(r')^2 - r r'' + r^2 \right).
\label{eq:curv}
\end{equation}
The right side of (\ref{eq:csre}) is just
\begin{equation}
F \cdot \nu = -\frac{r^2}{\sigma}.
\label{eq:rs}
\end{equation}
Equating these and simplifying yields, finally, the main equation
\begin{equation}
r r'' = r^2 + 2(r')^2 + r^2 (r^2 + (r')^2).
\label{eq:me}
\end{equation}

\begin{proposition}
Any maximally extended solution $r(\theta)$ of (\ref{eq:me}) is defined on an interval $(a,b) \subset 
[0,2\pi]$ (mod $2\pi$) with $b-a < \pi$, and satisfies:
\begin{itemize}
\item[i)] $r$ is convex as a function of $\theta$;
\item[ii)] $\lim_{\theta \searrow a}r(\theta) = \lim_{\theta \nearrow b} r(\theta) = \infty$;
\item[iii)] $r(\theta) = r((a+b) - \theta)$, i.e. the image of $r$ is symmetric about
the ray which makes angle $(a+b)/2$ with the horizontal.
\end{itemize}
\end{proposition}
\begin{figure}[ht]
\centering
\includegraphics[scale=0.5]{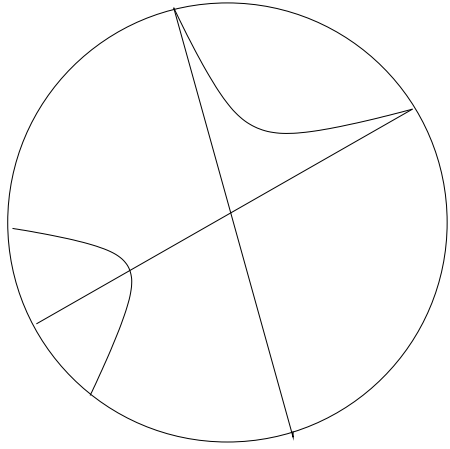}
\caption{Geodesics for $g$ in the stereographic compactification}
\end{figure}

\begin{proof} Property i) is obvious directly from (\ref{eq:me}). If $(a,b)$ is the maximal
interval of existence, then by convexity, the two limits in ii) exist. By earlier remarks,
neither limit can equal zero, and by maximality of $(a,b)$, these limits cannot be finite
either, which establishes ii). Using convexity and properness, we see that there is
a unique point $\theta_0 \in (a,b)$ where $r'(\theta_0) = 0$. Since (\ref{eq:me}) is
invariant under the reflection $\theta_0 + \eta \mapsto \theta_0 - \eta$, by uniqueness
of solutions of the initial value problem we see that the solution must be invariant
under this flip, and hence that $\theta_0 = (a+b)/2$. 

To prove that this maximal interval of existence has length less than $\pi$, define $v = r'/r$, so that 
\[
v' = \frac{r''}{r} - \frac{(r')^2}{r^2}.
\]
Using (\ref{eq:me}), we find that
\[
v' = (1+v^2)(1+r^2) \Longrightarrow v' \geq (1+v^2).
\]
Integrating gives
\[
\frac{r'(\theta)}{r(\theta)} \geq \tan (\theta + C)
\]
for some $C$, and this clearly proves the claim.

Note that if $r$ attains its minimum at $\theta_0$ and $r(\theta_0) = r_0$, then
\[
v' \geq (1+v^2)(1+r_0^2),
\]
hence
\[
v(\theta) \geq \tan ((1+r_0^2)\theta + C),
\]
so that the maximal interval of existence of this solution is of length less than $\pi/(1+r_0^2)$. 
\end{proof}

These solutions account for all remaining solutions of (\ref{eq:csre}). 

The most striking feature of these solutions is that they behave qualitatively exactly like
the geodesics on the hyperbolic plane. In particular, we have the

\begin{corollary}
Let $B$ be the compactification of the stereographic projection of $\RR^2$. Then every maximally 
extended solution of (\ref{eq:csre}) intersects $\del B$ in precisely two points. Moreover, for 
any pair of distinct points $p,q \in \overline{B}$, there exists a unique solution $\gamma$ 
of (\ref{eq:csre}) which passes through (or terminates) at these two points.
\end{corollary}

This behaviour is no accident, since in fact these solution curves are geodesics for a metric 
of negative curvature! 

\begin{proposition} 
Let $g_0 = dr^2 + r^2 d\theta^2$ be the standard Euclidean metric on $\RR^2$ and define $g = e^{r^2}g_0$. 
The geodesics for $g$ are solutions to (\ref{eq:csre}) and conversely.
\end{proposition}
\begin{proof}
Since the conformal factor $e^{r^2}$ is radial, the straight lines through the origin are geodesics for $g$. 
Suppose that $(r(u),\theta(u))$ is the polar representation of any $\calC^2$ curve in the plane. The 
geodesic equations for $g$ are
\begin{eqnarray*}
\ddot{r} - (r^3 + r) \dot{\theta}^2 + r \dot{r}^2 & = 0 \\
\ddot{\theta} + 2(r + r^{-1})\dot{r}\dot{\theta} & = 0
\end{eqnarray*}
(where now the dot refers to derivatives with respect to the parameter $u$). Assuming that $\dot{\theta} \neq 0$, we 
consider $r$ as a function of $\theta$ and immediately derive that
\[
r'' = r + 2r^{-1}(r')^2 + r(r^2 + (r')^2),
\]
which is just (\ref{eq:me}); here $r' = dr/d\theta$. 

The proof is complete.
\end{proof}

The Gauss curvature of the metric $g$ is equal to $K = -2e^{-r^2}$, which is everywhere negative.
However, since $K \to 0$ at infinity, we cannot circumvent the proposition above about the global 
behaviour of geodesics by appealing to the properties of geodesics on Cartan-Hadamard manifolds.

One final remark is that any two distinct geodesics of $g$ which converge to the same boundary
point $p \in \del B$ meet tangentially there, but the order of tangency is not quadratic.
Indeed, a short estimation shows that while the tangents to these curves do coincide at $p$, 
their directions approach one another inversely proportionally to the logarithm of the distance 
to the origin (with respect to $g$).

\section{Regular geodesic networks}
The remaining part of the proof of our main theorem involves showing the existence of a regular 
network in $\RR^2$ where each edge is a geodesic arc or ray for the metric $g$, and whose
tangent cone at infinity is $C_0$, the union of $k$ half-lines meeting at $0$. Alternately,
this network spans $k$ specified points $p_1, \ldots, p_k \in \del B$, where $B$ is the 
compactification of the stereographic projection of $\RR^2$. For convenience, we assume that 
these points are labelled in consecutive order around the circle (mod $k$). 

When $k=3$, the existence and uniqueness of this regular network may be accomplished directly by 
degree theory, cf.\ \cite{selsimos}. Certain special cases are also quite easy to handle: for 
example, if $k$ is even, then a particular solution is the disjoint collection of geodesics
$\gamma_1, \ldots, \gamma_{k/2}$,  where $\gamma_j$ connects $p_{2j-1}$ to $p_{2j}$, $j = 1, \ldots, k/2$. 
Such a network would correspond to a `complete dissolution' of this vertex into smooth curves.
We consider instead the other extreme where the solution is connected. The remainder
of this section is devoted to proving the existence of at least one connected regular geodesic
network which spans these $k$ points. The number of connected solutions is, by definition,
just the number of Steiner trees (with respect to $g$) with this given boundary, but the precise
number is a priori not so obvious.  We conjecture that this number is actually the same as
if we were taking Steiner trees with respect to the hyperbolic metric on the ball.

It is possible to construct connected regular geodesic networks with arbitrary prescribed boundary
values using synthetic geometry. However, another -- particularly efficient -- way to obtain existence 
uses geometric measure theory. The precise formulation from this point of view, as well as relevant 
literature, was suggested to us by Brian White, to whom we are very grateful. 

We shall work in the class $\calF_1(\RR^2;\ZZ_k)$ of flat chains of dimension $1$ in $\RR^2$ with coefficients 
in the group $\ZZ_k$, endowed with the norm $|g| = 1$ for all $g \in \ZZ_k$. Let us explain what this means. 
First, the space of flat chains of dimension $1$ is the completion of the space of polygonal 
curves with respect to the flat norm. A flat chain with coefficients in $\ZZ_k$ is an ordinary 
flat chain such that its (integer) multiplicity function is reduced mod $k$. The norm on $\ZZ_k$ 
appears in the definition of the {\it size } (rather than mass) of a flat chain $c = \sum g_\tau 
\tau$, where each $g_\tau \in \ZZ_k$, which is given by
\[
\SSS(c) = \sum_\tau |g_\tau| {\mathbb M}(\tau) = \sum_\tau {\mathbb M}(\tau).
\]
Here ${\mathbb M}(\tau)$ denotes mass with respect to the metric $g$, which for a $\calC^1$
arc corresponds to the usual length in that metric. We refer to \cite{W} and also \cite{Mo2} 
for these facts and for more complete references. 

For each $j$, let $\ell_j$ denote the ray from $0$ to the point $p_j$ at infinity, 
and for each radius $R$, set $p_j^R = \ell_j \cap \del B_R$. 

\begin{lemma}
There exists a connected $1$-dimensional flat chain $T^R$ with coefficients in $\ZZ_k$ such that 
$\del T^R=\{p_1^R, \ldots, p_k^R\}$ and which is size-minimizing, i.e.
\[
\SSS(T^R)=\inf\left\{\SSS(S): S \in \calF_1(B_R; \ZZ_k), \ \mbox{supp}(S)\, \mbox{connected},\ 
\del S= \{p_1, \ldots, p_k\} \right\}.
\]
The support of this minimizer $T^R$ is a regular geodesic network: each edge is a geodesic arc
or ray for the metric $g$ and there are a finite number of interior vertices, each of which is
trivalent, where these edges meet in equal angles. 
\end{lemma}

\begin{proof}
This is all standard, but we review the argument briefly. (There is one twist at the end,
about the valence of boundary vertices.) For the existence we use the compactness theorem for 
flat chains, as in \cite{W}, cf.\ also \cite{Mo1}. 
Let $T_j$ be a sequence of elements in $\calF_1(B_R;\ZZ_k)$ with connected support and $\del T_j 
= \{p_1, \ldots, p_k\}$, and such that $\SSS(T_j)$ converges to the minimum possible value in this 
class of competitors. Clearly $\SSS(T_j) \leq C$ for some $C>0$, and in addition $\SSS(\del T_j) = k$, 
for all $j$. The compactness theorem implies that there is a convergent 
subsequence, relabeled again as $T_j$, with limit $T^R$. This has the same boundary, and 
lower-semicontinuity of $\SSS$ implies that 
\[
\SSS(T^R)\leq \lim_{j\to \infty} \SSS(T_j),
\]
so that $T^R$ is indeed a size minimizer.

For the regularity of $T^R$, we use \cite{AA}, but see also \cite{Mo1}, \cite{Mo2} and \cite{T}. 
The regularity theorem in \cite{AA} implies that the singular set has Hausdorff dimension $0$, and in 
fact consists of a finite number of points. By the first variation formula, the regular set consists 
of a finite number of geodesic arcs or rays. At each interior vertex $p$, following \cite{Mo1}, 
consider the following variation. Let $v$, $w$ be unit vectors tangent to two adjacent edges 
$\ell_v$ and $\ell_w$ which 
meet at $p$. Let $u$ be another unit vector in the positive cone determined by $v$ and $w$, and 
consider the new network which replaces these two edges by a triod with one very short edge of length
$\epsilon$ along the geodesic starting at $p$ in the direction $u$ and the other two edges the 
geodesics from the other end of that short geodesic to the other ends of $\ell_v$ and $ \ell_w$. 
The multiplicity of the short geodesic should equal the sum of the multiplicities of $\ell_v$
and $\ell_w$, while the two new longer geodesics should have the same multiplicities as those,
respectively. We compute that
\[
0\leq \left.\frac{d\,}{d\epsilon} \SSS(T^R_\epsilon)\right|_{\epsilon=0}=1-(v+w)\cdot u;
\]
the inequality holds because $T^R$ is minimizing. In particular, setting $u=\frac{v+w}{|v+w|}$, 
we conclude that $1\leq |v+w|$, or equivalently $v\cdot w\leq-\frac{1}{2}$. Hence the angle 
between $v$ and $w$ is at least $2\pi/3$. Therefore, at most three edges can meet at $p$, 
and if there are three incoming vertices, then these must meet at $2\pi/3$. If only two 
edges meet at an interior vertex $p$, then the geodesic connecting the other two endpoints 
of these edges would be shorter. 

It can happen in certain geometries that precisely two vertices meet at a boundary point,
though of course by the preceding argument, the angle between them must be at least $2\pi/3$. 
Since $\del B_R$ is convex, neither of these edges will lie along this boundary. We claim that
this is impossible once $R$ is sufficiently large. Indeed, the convex hull of the points
$p_1^R, \ldots, p_k^R$ in $B_R$ is a polygon with geodesic sides, and with angle at
each vertex tending to $0$ as $R \to \infty$, cf.\ the final remark of \S 2. Hence these boundary 
vertices must be univalent as soon as $R$ is sufficiently large so that the opening angle of 
this convex hull is less than $2\pi/3$.

To conclude, we must show that $T^R$ is connected. Since convergence in flat norm implies convergence 
as currents (see \cite{S}), 
\begin{equation}
\int_{T_j}\phi \, dT_j\to \int_{T^R}\phi \, dT^R
\label{eq:convcurr}
\end{equation}
for every compactly supported $\phi$. Suppose that the support of $T^R$ is disconnected. Then 
there is a curve $\gamma$ in $B_R$ such that $B_R \setminus \gamma$ has two components,
each intersecting the support of $T^R$ nontrivially. Denote by $\calU_\e$ the $\e$-neighbourhood
around $\gamma$. For sufficiently small $\e$, $\calU_\e \cap \mbox{supp}\,(T^R) = \emptyset$.
Consider a nonnegative $\phi \in \calC^\infty_0$ with support in $\calU_\e$ which equals $1$
in $\calU_{\e/2}$. Since $T_j$ is connected  
\[
\int_{T_j}\phi \, dT_j\geq \frac{\e}{4}
\]
for all $j$, but on the other hand 
\[
\int_{T^R} \phi \, dT^R = 0,
\]
which contradicts (\ref{eq:convcurr}). This finishes the proof.
\end{proof}

To obtain a network which spans the points $p_1, \ldots, p_k \in \del B$, we take the limit of
$T^R$ as $R \to \infty$. 

\begin{proposition}
Let $R_j$ be a sequence of radii tending to infinity, and let $T_j$ be one of the connected size-minizing 
flat chains with coefficients in $\ZZ_k$ with $\del T_j = \{p_1^{R_j}, \ldots, p_k^{R_j}\}$ obtained in the 
previous lemma. As $j \to \infty$, some subsequence of the $T_j$ converges (in the flat topology) to a 
locally size-minimizing connected flat chain with $\del T = \{p_1, \ldots, p_k\}$.
\end{proposition}
\begin{proof}
Let $\calP$ denote the ideal $k$-gon which is the convex hull of the points $p_1, \ldots, p_k$,
and $\calP_j$ the convex hull of $p_1^{R_j}, \ldots, p_k^{R_j}$. As already used in the last proof, 
the support of each $T_j$ lies in $\calP_j \subset \calP$, and $\calP_j \nearrow \calP$. 

Next, let us observe that the total number of interior vertices in the support of each $T_j$
remains fixed. Indeed, if $\ell$ denotes the number of interior vertices and $e$ the number of edges, then
\begin{eqnarray*}
3\ell + k & = 2e \\
\ell + k & = e + 1.
\end{eqnarray*}
The first equation uses that each interior vertex is trivalent and each boundary vertex connects to
only one edge; the second equation asserts that the Euler characteristic of a tree is equal to $1$.
Subtracting the second equation from the first gives $2\ell = e-1$, and hence, after some manipulation
\[
\ell = k-2.
\]

We now claim that no interior vertex can converge to any one of the boundary vertices, and
hence disappear in the limit. Indeed, if an interior vertex $q$ lies in the cusp of this convex
hull corresponding to some $p_i$, then at most one of the three edges which meet at $q$ is directed
`outwards', toward $p_i$ and the others must be pointed inward. However, since the tree remains
inside the convex hull, the extensions of either of the other edges must hit the boundary of the
convex hull in some distance which we can estimate from the position of $q$. Since this is impossible,
there must be two new interior vertices. We can repeat this argument a finite number of times. If
$q$ is sufficiently deep into this corner, the tree would have more than $k-2$ additional vertices
in just this neighbourhood, which is impossible. 

This shows that all $k-2$ interior vertices remain within some fixed ball $B_{R}$, and hence
we can take a limit of the $T_j$ and obtain a nontrivial limit $T$.  Clearly $\del T =
\{p_1, \ldots, p_k\}$, as required, and by the same argument as above, the support of $T$
is connected. Finally, it is a standard fact that the limit of a convergent sequence of 
mass-minimizing currents is again mass-minimizing. This transfers immediately to the setting
of flat chains with coefficients in $\ZZ_k$, with the given norm on this cyclic group. 
\end{proof}

It is worth remarking that this result (and argument) is close in spirit to the construction 
in \cite{An} of complete mass-minimizing submanifolds in hyperbolic space which have a prescribed 
asymptotic boundary at infinity. 

\section{Behaviour of the flow at $t=0$}
We conclude this paper with some brief remarks about the precise regularity of the self-similar solutions
to this network flow which we have constructed here. This will be expanded on considerably
in our subsequent paper on general short-term existence results, where it plays a more crucial role. 

Even when $k=3$ but the initial half-lines from $0$ do not meet at equal angles, there seems to be 
a suddent jump in the configuration as soon as $t$ becomes positive. Of course, when $k > 3$, this 
jump is even more pronounced since new vertices and edges are created instantaneously.  There is a 
way of viewing all of this, however, which makes this behaviour continuous. To this end, we introduce 
the parabolic blowup of $\RR^2 \times \RR^+$ at $\{x=y=t=0\}$. This is a manifold with corners of 
codimension two which is obtained by taking the union of $\{(x,y,t): t \geq 0, (x,y,t) \neq (0,0,0)\}$ 
and a new hypersurface $F$, which is the `parabolic spherical normal bundle' of the origin with respect to the 
family of dilations $D_\lambda$ introduced in \S 2. In other words, every point of $F$ corresponds to 
an orbit of this dilation group. If we were using ordinary dilations ($(x,y,t) \to (\lambda x, 
\lambda y, \lambda t)$) then this would be the more familiar normal blowup, which can
be described easily in terms of polar coordinates around the origin: indeed, the new face
added in a normal blowup is the one obtained by setting the radial variable equal to $0$.
The picture is the same, however, since the new face is diffeomorphic to a half-sphere.
We denote this blown up space by $X$; the two codimension one boundaries are the new
face $F$, described above, and the compactification of the original boundary minus
the origin, which we denote $T$. Note that $T$ is naturally the complement of
a ball in $\RR^2$, while $F$ is a half-sphere.

In the previous sections, we were using an identification of the slice $\{t=1\}$ with an open 
hemisphere via ordinary stereographic projection. There is a similar identification of this 
slice with the interior of the face $F$ defined by the dilations $D_\lambda$. Note, however,
that this identification is {\it not} conformal. 

The world-sheet $\Gamma$ of any self-similar solution of the network flow is a union
of pieces of surfaces in $\RR^2 \times \RR^+$, where these smooth components intersect
along curves which are orbits of the family of dilations. The entire surface $\Gamma$
is a union of such dilation-invariant curves, and hence the closure of 
$\Gamma \setminus (0,0,0)$ in $X$ intersects $F$ and $T$ in a certain collection of
curves. The intersection with $T$ is just the union of half-lines, while the intersection
with $F$ is a regular network on the hemisphere. The intersections of these curves
at the corner $F \cap T$ yield the $k$ boundary points $p_1, \ldots, p_k$. We illustrate this 
when $k=4$. 

\begin{figure}[ht]
\centering
\includegraphics[scale=0.5]{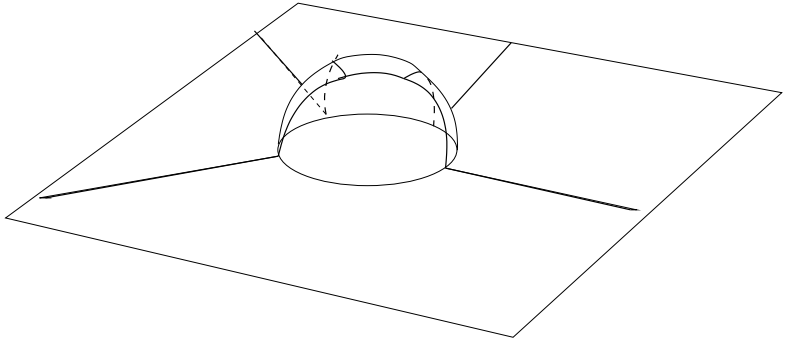}
\caption{The blowup and induced network at $t=0$}
\end{figure}

Our final result is tautological, but is an important key in understanding regularity
near $t=0$ for the network flow with more general initial configurations.
\begin{proposition}
The closure of $\Gamma \setminus \{(0,0,0)\}$ in $X$ is a union of a finite number
of smooth surfaces, each with boundary and corners, which intersect along a finite
number of smooth curves. Each of these curves of intersection is an orbit of the 
family of dilations. 
\end{proposition}


\begin{thebibliography}{MNT04}
\bibitem[AA]{AA}
W. Allard and F. J. Almgren, Jr.
\newblock The structure of stationary one dimensional varifolds with positive density
\newblock {\em Invent. Math.} {\bf 34} (1976), 83-97. 

\bibitem[An]{An}
M. Anderson
\newblock Complete minimal varieties in hyperbolic space
\newblock {Invent. Math.} {\bf 69} (1982), 477-494.

\bibitem[ACI]{ACI}
S.B. Angenent, D. Chopp, and T. Ilmanen
\newblock  A computed example of nonuniqueness of mean curvature flow 
in $\RR^3$
\newblock Comm. Part. Diff. Eq. {\bf 20} (1995) 1937-1958.

\bibitem[Br]{Br}
K. Brakke
\newblock The motion of a surface by its mean curvature
\newblock Mathematical Notes \# 20, Princeton University Press, Princeton (1978).
 
\bibitem[HS]{HS}
N. Hungerb\"uhler and K. Smoczyk 
\newblock {\em Soliton solutions for the mean curvature
flow}, Diff. and Int. Eqns. {\bf 13} No. 10-12 (2000), 1321-1345.

\bibitem[I]{I}
T. Ilmanen
\newblock Lectures on mean curvature flows and related equations
\newblock http://www.math.ethz.ch/~ilmanen/papers/notes.pdf

\bibitem[MNT]{moybycm}
C. Mantegazza, M. Novaga, and V.~M. Tortorelli.
\newblock Motion by curvature of planar networks.
\newblock {\em Ann. Sc. Norm. Super. Pisa Cl. Sci. (5)} {\bf 3} No. 2 (2004), 235--324.

\bibitem[M1]{Mo1}
F. Morgan
\newblock Geometric Measure Theory: A beginner's guide.
\newblock {\em New York: Academic Press, 1988}
 
\bibitem[M2]{Mo2}
F. Morgan
\newblock Size-minimizing rectifiable currents.
\newblock {\em Inv. Math.} {\bf 96} (1989), 333--358.

\bibitem[M3]{Mo3}
F. Morgan
\newblock Soap Bubbles in $\RR^2$ and in surfaces.
\newblock {\em Pacific Jour. Math.} {\bf 165} No. 2 (1994), 347--361.

\bibitem[SS]{selsimos}
O.~Schn\"urer and F.~Schulze.
\newblock Self-similar expanding networks to curve shortening flow.
\newblock arXiv:math-DG/0702698

\bibitem[S]{S}
L. Simon.
\newblock Lectures on geometric measure theory.
\newblock Australian National University, Canberra (1983).
 
\bibitem[T]{T}
J. Taylor.
\newblock The structure of singularities in soap-bubble-like and soap-film-like minimal surfaces.
\newblock {\em Ann. of Math.} {\bf 103} (1976), 489--539

\bibitem[W]{W}
B. White
\newblock Existence of least-energy configurations of immiscible fluids
\newblock {\em Jour. Geom. Anal.} {\bf 6} No. 1 (1996), 151-161.
\end{thebibliography}
\end{document}